\newtheorem{defn0}{Definition}[section]
\newtheorem{prop0}[defn0]{Proposition}
\newtheorem{thm0}[defn0]{Theorem}
\newtheorem{lemma0}[defn0]{Lemma}
\newtheorem{corollary0}[defn0]{Corollary}
\newtheorem{example0}[defn0]{Example}

\newenvironment{prop}{\begin{prop0}}{\end{prop0}}
\newenvironment{thm}{\begin{thm0}}{\end{thm0}}
\newenvironment{lemma}{\begin{lemma0}}{\end{lemma0}}
\newenvironment{corollary}{\begin{corollary0}}{\end{corollary0}}
\newenvironment{example}{\begin{example0} \rm }{\end{example0}}

\newcommand{\propref}[1]{Proposition~\ref{#1}}

\newcommand{\lemref}[1]{Lemma~\ref{#1}}
\newcommand{\corref}[1]{Corollary~\ref{#1}}
\newcommand{\exref}[1]{Example~\ref{#1}}

\newcommand{\qed}{\mbox{~~~~\vrule height 1.2ex width .9ex depth .1ex}}
\newenvironment{proof}{\noindent {\bf Proof.}}{\qed\vskip 6pt}

\newcommand{\std}{Gr\"{o}bner}
\newcommand{\Std}{Gr\"{o}bner}

\newcommand{\mmbox}[1]{\mbox{${#1}$}}

\newcommand{\Spec}[1]{\mmbox{{\rm Spec\,}{#1}}}

\newcommand{\iin}[1]{\mmbox{{\rm in}({#1})}}
\newcommand{\Ann}[1]{\mmbox{{\rm Ann}({#1})}}
\newcommand{\caps}[3]{\mmbox{{#1}_{#2} \cap \ldots \cap {#1}_{#3}}}

\newcommand{\sliver}{\hskip 0.01in}
\newcommand{\ab}{\sliver{\bf a}\sliver}
\newcommand{\x}{\sliver{\bf x}\sliver}

\newcommand{\Ax}{\mmbox{A[\x]}}
\newcommand{\Bx}{\mmbox{B[\x]}}
\newcommand{\kpx}{\mmbox{k(p)[\x]}}
\newcommand{\Apx}{\mmbox{A_p[\x]}}
\newcommand{\N}{{\bf N}}
\newcommand{\Oh}{{\cal O}}
\newcommand{\Z}{{\bf Z}}

\newcommand{\quotes}[1]{\mbox{``#1''}}

\documentstyle[11pt]{article}
\title {\Std\ Bases and Extension of Scalars}

\author
{
Dave Bayer
\thanks {Partially supported by
  the Alfred P. Sloan Foundation,
  by ONR contract N00014-87-K0214,
  and by NSF grant DMS-90-06116.}
\and Andr\'e Galligo
\and Mike Stillman
\thanks {Partially supported by
  the U.S. Army Research Office through
  the Mathematical Sciences Institute of Cornell University,
  and by NSF grants CCR-89-01061 and DMS-88-02276.}
}

\date{September 30, 1991}

\begin{document}

\maketitle

\begin{center}
{\em
 Dedicated to Professor Heisuke Hironaka on his sixtieth
birthday.
}
\end{center}

\section{Introduction}
\label{intro}
\setcounter{defn0}{0}

Let $A$ be a Noetherian commutative ring with identity, let $\Ax =
A[x_1,\ldots ,x_n]$ be a polynomial ring over $A$, and let $I \subset
\Ax$ be an ideal. Geometrically, $I$ defines a family of schemes over
the base scheme \Spec{A}; the fiber over each point $p \in \Spec{A}$
is a subscheme of the affine space ${\bf A}^n_{k(p)} =
\Spec{k(p)[\x]}$, where $k(p) = A_p/p_p$ is the residue field of $p$.

Let $>$ be a total order on the monomials of \Ax\ satisfying $\x^E >
\x^F \Rightarrow \x^G\x^E > \x^G\x^F$, and satisfying $x_i > 1$ for
each $i$. For $f \in \Ax$, define \iin{f} to be the initial (greatest)
term $c \x^E$ of $f$ with respect to the order $>$, where $c \in A$ is
nonzero. For $I \subset \Ax$, define the initial ideal $\iin{I}$ to be
the ideal $( \iin{f} \mid f \in I )$ generated by all initial terms of
elements of $I$. $\iin{I} \subset \Ax$ is generated by single terms of
the form $c \x^E$; we call such an ideal a monomial ideal.
$\{f_1,\ldots , f_r\} \subset I$ is a \std\ basis for $I$ if and only
if $\{\iin{f_1},\ldots , \iin{f_r}\}$ generates \iin{I}.

For an ideal $I$ in a ring of formal power series $A[[\x]]$, Hironaka
defined the corresponding notion (the standard basis of $I$) and
proved a generalized Weierstrass division theorem.The relation between
flatness and division was considered in \cite{hlt73} and later in
\cite{gal79} in order to obtain a presentation for the flattener of a
germ of an analytical morphism. The algebraic situation is slightly
different.

In this paper, we study the behavior of \std\ bases with respect to an
extension of scalars $A \rightarrow B$. When does a \std\ basis for
$I$ map to a \std\ basis for $I \Bx$? It suffices to have \iin{I}
generate \iin{I \Bx}; we focus on this condition. Taking $B = k(p)$
for $p \in \Spec{A}$, we consider the relationship between a \std\
basis for $I$, and the \std\ bases of the fibers of the family defined
by $I$. How much information about the fibers of this family can be
inferred from knowledge of \iin{I} alone?

Let $X \subset \Spec{A}$ be the support of the family defined by $I$.
A \Std\ basis for $I$ encodes considerable information about this
family, even when $X$ is nonreduced or reducible. To interpret \iin{I}
in such situations, we work with its {\it coefficient ideals\,}: The
coefficient ideal for the monomial $\x^E$ vanishes on the support of
those fibers where $\x^E$ fails to belong to $\iin{I} \kpx$. From this
point of view, a point $p \in \Spec{A}$ is \quotes{good} if each
coefficient ideal of \iin{I} defines a scheme which either avoids $p$,
or contains an open neighborhood of $p$ in $X$.

In \S 2, we study coefficient ideals of monomial ideals. In \S 3, we
prove that an extension of scalars commutes with taking the initial
ideal of any ideal $I$, if and only if the extension is flat. We then
prove that a \std\ basis for $I$ determines \std\ bases for the
localizations to dense open subsets of each isolated component of $X$.
We also prove that for this family, \iin{I} determines the fiber
initial ideals over \quotes{good} points, as defined above. These
results reveal that \iin{I} carries generic information for each
isolated component of $X$. In \S 4, we prove that if every point of
$X$ is \quotes{good}, then the family defined by $I$ is faithfully
flat over $X$. Faithful flatness imposes strong conditions on the
component structure of the total space of our family, so this result
has geometric applications, such as the removal of unwanted
components. In \S 5, we give two other applications of coefficient
ideals, describing the locus where a morphism of schemes is an
isomorphism, or a finite map.

If $A$ is a finitely generated $k$-algebra for a field $k$, then we
can write $A = k[a_1,\ldots,a_m]/J$ for some ideal $J$. We can
reformulate our problem as $A = k[\ab]$, with $I \subset \Ax$ and $I
\cap A \supset J$. A \std\ basis for $I$ can then be computed by the
usual algorithm over a field, by combining orders $>_1$, $>_2$ into a
product order $$ \ab^D \x^E > \ab^F \x^G \ \Leftrightarrow \x^E >_1
\x^G, \mbox{\ or \ } \x^E = \x^G \mbox{\ and \ } \ab^D >_2 \ab^F.$$ In
this setting, $I$ defines a subscheme $Y \subset {\bf A}_k^{m+n}$
which projects to $X \subset {\bf A}_k^m$.

More generally, the computational relevance of this work depends on
our ability to compute in the base ring $A$. Specifically, $A$ needs
to be a ring where linear equations are solvable; see Trager, Gianni,
and Zacharias (\cite{gtz88}) for background material and references on
\std\ bases in this setting. Our paper continues their study of
families of \std\ bases; we would like to thank each of them, David
Eisenbud, and an anonymous referee, for many helpful conversations and
suggestions.


\section{Monomial Ideals}
\label{monoms}
\setcounter{defn0}{0}

Let $J \subset \Ax$ be a monomial ideal, i.e., an ideal generated by
single terms of the form $c \x^E$, with $c \in A$.

When $A$ is a field, $J$ is easily understood: its structure is
realized by the subset $L = \{ E \mid \x^E \in J \}$ of $\N^n$, where
$\N$ denotes the nonnegative integers. $\N^n$ admits a natural partial
order $\leq$ defined by $E \leq F$ iff $\x^E$ divides $\x^F$. The
characteristic function of $L$ can be viewed as a poset homomorphism
from $\N^n$, ordered by $\leq$, to the set of ideals $\{(0),(1)\}$ of
$A$, ordered by inclusion.

To understand $J$ when $A$ is not a field, it is helpful to consider
the {\it coefficient ideals\,} of $J$: Define $J_E = J_{\x^E} = (c \in
A \mid c \x^E \in J)$. Alternatively, $J_E$ is the ideal quotient $(J
: \x^E) \cap A$. This construction defines a poset homomorphism $E
\mapsto J_E$ from $\N^n$, ordered by $\leq$, to the set of ideals of
$A$, ordered by inclusion. Conversely, any such poset homomorphism
determines a monomial ideal $J$, so we can think of this construction
as describing an equivalence of categories.

Viewing a monomial ideal as its collection of coefficient ideals is
thus a purely tautological change of perspective; any operation on $A$
can be viewed as acting on $J$ via the inclusion $J \subset \Ax$, or
equivalently as acting on the set of coefficient ideals of $J$. In
particular, if $v: A \rightarrow B$ is a ring homomorphism, then $v$
extends naturally to a homomorphism $v: \Ax \rightarrow \Bx$. The
image under $v$ of any ideal $I \subset \Ax$ generates the extension
ideal $I^e = I \Bx$. For monomial ideals, we have the immediate
proposition

\begin{prop}
\label{monoms:prop}
Let $v: A \rightarrow B$ be a ring homomorphism. Let $J B[\x]$ denote
the monomial ideal obtained by extension of scalars from $A$ to $B$,
and let $J_E B$ also be obtained by extension of scalars, for an
exponent $E$. Then $$J_E B = (J \Bx)_E.$$
\qed
\end{prop}

In particular, if $B$ is the residue field $k(p)$ of a prime ideal $p
\subset A$, then \propref{monoms:prop} asserts that $\x^E \in J \kpx$
if and only if the point $p$ does not belong to the subscheme of
\Spec{A} defined by $J_E$. This subscheme is the support of the
$A$-module $A \x^E \subset \Ax/J$. The monomial $\x^E$ appears only
with a coefficient of zero in $J (A/J_E)[\x]$; $J_E$ is the
intersection of all ideals $K \subset A$ with the property that $(J
(A/K)[\x])_E = (0)$.

In other words, we can view $J$ as a family of monomial ideals over
\Spec{A}. The monomial ideals corresponding to each fiber of the
family defined by $J$ are defined over fields, and can be visualized
combinatorially: Each monomial $\x^E$ either belongs or does not
belong to a given fiber monomial ideal, which in turn is determined by
this data. The coefficient ideal $J_E$ defines the subscheme of
\Spec{A}\ over which $\x^E$ does not belong to the fiber monomial
ideals. $J$ is determined by these subschemes.

\begin{example}
\label{monoms:ex1}
Let $A = \Z$, let $\Ax = A[x,y]$, and let $J = (9x,\, 2y,\, x^2,\,
y^2)$. The coefficient ideals for $1$, $x$, $y$, $x^2$, $xy$, and
$y^2$ respectively are $(0)$, $(9)$, $(2)$, $(1)$, $(1)$, and $(1)$;
this is summarized in the diagram

\begin{center}
\vspace{.1 in}
\begin{tabular}{r|lll}
$y$ & $(1)$ \\
    & $(2)$ & $(1)$ \\
    & $(0)$ & $(9)$ & $(1)$ \\
\cline{2-4}
\multicolumn{4}{r}{$x$}
\end{tabular}
\end{center}

$J$ specializes to $(x,\, y)$ in each fiber over the open subset of
\Spec{\Z} which is the complement of the points $(2)$ and $(3)$. In
the fiber over $(2)$, $J$ specializes to $(x,\, y^2)$. Over the double
point at $(3)$, $J$ specializes to $(x^2,\, y)$.

$A$ is an integral domain. The union of the monomials with a nonzero
coefficient ideal spans the monomial ideal $(x,\, y)$, which occurs
generically. On the other hand, the union of the monomials with
coefficient ideal $(1)$ spans the monomial ideal $(x^2,\, xy,\, y^2)$.
There is no specialization which produces this monomial ideal, but it
is contained in every monomial ideal obtained by specialization.
\end{example}

\begin{example}
\label{monoms:ex2}
Modify the preceding example by taking $A = \Z/18\Z$. \Spec{A} is no
longer reduced or irreducible. The union of the monomials with nonzero
coefficient ideals again spans the monomial ideal $(x,\, y)$. There is
no specialization to a field which produces this monomial ideal; this
can happen whenever $A$ is not an integral domain.
\end{example}


\section{Initial Ideals}
\label{initial}
\setcounter{defn0}{0}

Let $I \subset \Ax$ be an ideal. For monomial ideals,
\propref{monoms:prop} asserts that the formation of coefficient ideals
commutes with extension of scalars. In contrast, the formation of the
initial ideal of an arbitrary ideal $I$ need not commute with
extension of scalars: If $v: A \rightarrow B$ is a ring homomorphism,
then it can happen that $\iin{I} \Bx \neq \iin{I \Bx}$. This is
because if $v$ maps the leading coefficient of a polynomial $f \in I$
to zero, then the first surviving term of the image of $f$ will
contribute to \iin{I \Bx}. When this happens, \iin{I \Bx} cannot be
predicted from knowledge of \iin{I} alone.

Let $\{ f_1,\ldots f_r \}$ be a \std\ basis for $I$. If $\iin{I} \Bx =
\iin{I \Bx}$, then the images $\{ v(f_1),\ldots v(f_r) \}$ form a
\std\ basis for $I$. We shall study the behavior of initial ideals
with respect to various hypotheses on $I$ and $B$; as a consequence we
obtain sufficient conditions for the construction of a \std\ basis for
$I$ to commute with the extension of scalars $v: A \rightarrow B$.

\begin{example}
\label{initial:ex1}
Let $A=k[a]$ for a field $k$, let $\Ax=A[x,y]$, and let $I = (ax-y)
\subset \Ax$. In this example, as in all subsequent examples involving
the variables $x$ and $y$, we use the lexicographic order extending $x
> y$. Choose a prime $p \subset A$, and let $k(p)$ be the residue
field $A_p/p_p$ of $p$. When $p \neq (a)$, $\iin{I} \kpx = \iin{I
\kpx} = (x)$. However, when $p = (a)$, $\iin{I} \kpx = (0)$, but
$\iin{I \kpx} = (y)$. Nevertheless, the image of $\{ ax-y \}$ is a
\std\ basis for $I \kpx$.

$I$ defines a faithfully flat family over \Spec{A}, because $\Ax/I
\simeq A[x]$. The total space is given by the surface $a x - y = 0$,
and each fiber consists of a line through the origin in ${\bf A}^2_k$,
with slope $a$. \iin{I \kpx} momentarily flips from $(x)$ to $(y)$ as
this slope passes through zero.
\end{example}

\begin{example}
\label{initial:ex2}
Let $A = k[a,b]$, let $\Ax = A[x,y]$, let $I = (ax^2+y,\, by^2+y+1 )$,
and let $B = A/(a,\, b)$. Then $\iin{I} \Bx = (0)$, but $\iin{I \Bx} =
(1)$. $\{ax^2+y,\, by^2+y+1\}$ is a \std\ basis for $I$, but its image
$\{y,\, y+1\}$ in $\Bx$ is not a \std\ basis for $I \Bx$.
\end{example}

\begin{example}
\label{initial:ex3}
Let $A = k[a]$, let $\Ax = A[x]$, let $I = (ax - 1)$, and let $p
\subset A$ be a prime. When $p \neq (a)$, $\iin{I} \kpx = \iin{I \kpx}
= (x)$. However, when $p = (a)$, $\iin{I} \kpx = (0)$, but $\iin{I
\kpx} = (1)$.

$I$ defines a flat family which is not faithfully flat: The fiber over
$a = 0$ is empty, so $(a) \subset A$ extends to the unit ideal in
$\Ax/I$.
\end{example}

As suggested by these examples, we do have an inclusion in one direction:

\begin{prop}
\label{initial:incprop}
For any ring homomorphism $v: A \rightarrow B$, and for any ideal $I
\subset \Ax$, we have $$\iin{I} \Bx \subset \iin{I \Bx}.$$
\end{prop}

\begin{proof}
It is enough to show that each generator of $\iin{I} \Bx$ also belongs
to $\iin{I \Bx}$. $\iin{I} \Bx$ is generated by $v(\iin{f})$ for $f
\in I$. For each $f \in I$, either \iin{f} maps to zero in \Bx, or
else $v(\iin{f}) = \iin{v(f)} \in \iin{I \Bx}$. \end{proof}

The following theorem asserts that taking initial ideals universally
commutes with an extension of scalars if and only if the extension is
flat. We apply the criterion for flatness given in \cite{mat86}, Thm.\
7.6, which asserts that $v: A \rightarrow B$ is flat iff the syzygies
in $B$ of a set of elements from $A$ can always be generated by
syzygies from $A$:

\begin{lemma}
\label{initial:flatlem}
Let $v: A \rightarrow B$ be a ring homomorphism. $v$ is flat if and
only if for each sequence $a_i \in A$ and $b_i \in B$ for $1 \leq i
\leq r$ so $$\sum_i b_i v(a_i) = 0,$$ then for some $s$ we can choose
$c_{ij} \in A$ and $d_j \in B$ for $1 \leq j \leq s$ so $$\sum_i
c_{ij} a_i = 0 \mbox{\ for each $j$, and\ } b_i = \sum_j d_j v(c_{ij})
\mbox{\ for each $i$.}$$
\qed
\end{lemma}

\begin{thm}
\label{initial:thm}
Let $v: A \rightarrow B$ be a ring homomorphism. The following two
conditions are equivalent:

(a) for any number of variables $x_1, \ldots, x_n$, and for any ideal
$I \subset \Ax$, $\iin{I} \Bx = \iin{I \Bx}$;

(b) $B$ is a flat $A$-algebra.
\end{thm}

\begin{proof}
First, suppose that (b) holds, and let $I \subset \Ax$ be an ideal. We
need to show that $\iin{I \Bx} \subset \iin{I} \Bx$. Given $c \x^E \in
\iin{I \Bx}$, consider expressions of the form $c \x^E = \iin{\sum_i
b_i v(f_i)}$, where $b_i \in \Bx$ and $f_i \in I$. By expanding out
each $b_i$ and absorbing variables into each $f_i$, we need only
consider expressions for which each $b_i \in B$. Among all such
expressions, choose one for which the greatest monomial appearing in
any of the $f_i$ is minimal. We claim that this greatest monomial is
$\x^E$. Letting $c_i$ be the coefficient of $\x^E$ in each $f_i$, we
have $c = \sum_i b_i v(c_i)$, so $c \x^E \in \iin{I}\Bx$.

Suppose otherwise, that the greatest monomial of a minimal expression
is $\x^D > \x^E$. Let $a_i$ be the coefficient of $\x^D$ in each
$f_i$. Then $\sum_i b_i v(a_i) = 0$. Choosing $c_{ij}$, $d_j$ as in
\lemref{initial:flatlem}, define $g_j = \sum_i c_{ij} f_i \in I$ for
$1 \leq j \leq s$. Then $$\sum_j d_j v(g_j) = \sum_i b_i v(f_i),$$ and
$\iin{g_j} < \x^D$ for each $j$, contradicting the minimality of our
expression. This proves (a).

Now, suppose instead that (b) does not hold. Choose a sequence $a_i
\in A$ and $b_i \in B$ for $1 \leq i \leq r$ with $r$ minimal, so
$\sum_i b_i v(a_i) = 0$ but the $b_i$ cannot be expressed as in
\lemref{initial:flatlem}. We construct an example in two variables $x
> y$ for which $v$ does not commute with taking initial ideals: Let
$$f_i = a_i x^r + x^{r-i}y^i \mbox{\ \ for \ } 1 \leq i \leq r,$$ and
let $I \subset A[x,y]$ be defined by $ I = (f_1,\, \ldots,\, f_r)$.
Then $$\sum_i b_i v(f_i) = b_1 x^{r-1}y + \ldots + b_r y^r.$$
Moreover, for $c_{ij} \in A$, $$\sum_i c_{ij} a_i = 0
\;\Leftrightarrow\; c_{1j} x^{r-1}y + \ldots + c_{rj} y^r = \sum_i
c_{ij} f_i\in I.$$ Because the $b_i$ cannot be expressed as in
\lemref{initial:flatlem}, and because $r$ was chosen to be minimal,
already $b_1$ is not in the ideal generated by the images of all
$c_{1j}$ for $c_{ij}$ as above. Because $x^r > x^{r-1}y > \ldots >
y^r$, this ideal is the coefficient ideal of $\iin{I} \Bx$ with
respect to the monomial $x^{r-1}y$, so $$(\iin{I} \Bx)_{x^{r-1}y} \neq
\iin{I \Bx}_{x^{r-1}y}.$$ This proves that (a) does not hold.
\end{proof}

The following corollary asserts that taking initial ideals commutes
with taking rings of fractions, and is due to Gianni, Trager, and
Zacharias (\cite{gtz88}, Prop.\ 3.4).

\begin{corollary}
\label{initial:gtzcor}
Let $I \subset \Ax$ be an ideal, and let $B = S^{-1}A$ for a
multiplicatively closed set $S \subset A$. Then $$ \iin{I} \Bx =
\iin{I \Bx}.$$
\end{corollary}

\begin{proof}
A ring of fractions is a flat extension (\cite{mat86}, Thm.\ 4.5).
\end{proof}

Suppose that we want to work with the extension $A/J$, for an ideal $J
\subset A$. If $J$ arises as the kernel of a map $A \rightarrow
S^{-1}A$ for some multiplicatively closed set $S \subset A$, then we
can apply \corref{initial:gtzcor} if we instead work with the
extension $S^{-1}A$. Viewing $S^{-1}A$ as a ring of fractions of
$A/J$, this extension retains generic information along the scheme
defined by $J$, but loses primes annihilated by elements of $S$. Such
primes can prevent the taking of initial ideals from commuting with
the extension $A/J$, as illustrated by the following example.

\begin{example}
\label{initial:fractex}
Let $A = k[a,b]/(ab)$, let $\Ax = A[x]$, let $I = (ax+1)$, and let $J
= (a)$. Then $\iin{I} = (ax,\, b)$. Taking $B = A/J$, we have $\iin{I}
\Bx = (b)$, and $\iin{I \Bx} = (1)$. Instead taking $B = A_b$, we have
$\iin{I} \Bx = \iin{I \Bx} = (1)$.  Thus taking initial ideals
commutes with the extension $A_b$, but does not commute with the
extension $A/J$.

 \Spec{A_b} differs from \Spec{A/J} only by the removal of the prime
$p = (b)$. This prime obstructs good behavior for the extension $A/J$:
$\iin{I} \kpx = (0)$, and $\iin{I \kpx} = (1)$.
\end{example}

Which kernels arise from taking rings of fractions? From the proof of
\cite{atma69}, Thm.\ 4.10, one sees that these kernels are precisely
the ideals $q \subset A$ which arise as the intersection of primary
components corresponding to an isolated set of associated primes of
$(0)$. For each such $q$, it is enough to choose a multiplicatively
closed subset $S \subset A$ which intersects \Ann{q}, for $B =
S^{-1}A/S^{-1}q$ to equal $S^{-1}A$.

More generally, we may wish to consider components of the subscheme
defined by $I \cap A$, when this ideal is nonzero. The following lemma
reduces us to the above setting.

\begin{lemma}
\label{initial:quolem}
Let $I \subset \Ax$ be an ideal, and let $B = A/(I \cap A)$. Then $$
\iin{I} \Bx = \iin{I \Bx}.$$
\end{lemma}

\begin{proof}
Let $v:A \rightarrow B$ be the quotient map. We need to show that
$\iin{I \Bx} \subset \iin{I} \Bx$. Given $c \x^E \in \iin{I \Bx}$,
choose $f \in I \Bx$ so $\iin{f} = c \x^E$. Among all $g \in I$ so
$v(g) = f $, choose one with miminal leading term \iin{g}. We claim
that $v(\iin{g}) = \iin{f}$, so $c \x^E \in \iin{I} \Bx$.

Suppose otherwise, that $\iin{g} = b \x^D$ with $\x^D > \x^E$, and
$v(b) = 0$. Then $b \in I \cap A$, so $b \x^D \in I$, and $g - b \x^D
\in I$ has image $f$. This element has a lower leading term than $g$,
contradicting the minimality of our choice for $g$.
\end{proof}

\begin{prop}
\label{initial:genericprop}
Let $I \subset \Ax$ be an ideal, let the ideal $q \subset A$ be an
isolated primary component of $I \cap A$, and let the ideal $p \subset
A$ be its associated minimal prime. Define $B = A_p/q_p$. Then
$$\iin{I} \Bx = \iin{I \Bx}.$$
\end{prop}

\begin{proof}
Let $q \cap \caps q2s$ be a minimal primary decomposition of $I \cap
A$, with associated primes $p$, $p_2, \ldots, p_s$. Then $\caps q2s
\not\subset p$, for otherwise we would have $q_i \subset p$ for some
$i$, and thus $p_i \subset p$, contradicting the minimality of $p$.
Choose an element $r \not\in p$ such that $r \in \caps q2s$. Then $rq
\subset q \cap \caps q2s = I \cap A$, so $r \in (I \cap A:q)$. Thus
$A_p/q_p = A_p/(I \cap A)_p$.

By \lemref{initial:quolem}, taking initial ideals commutes with taking
the quotient by $I \cap A$. By \corref{initial:gtzcor}, taking initial
ideals commutes with forming a ring of fractions. The proposition
follows by combining these results.
\end{proof}

 \propref{initial:genericprop} affirms the utility of \std\ bases when
\Spec{A} is reducible: Enough information is encoded in such a basis
to determine the corresponding \std\ bases over dense open subsets of
each isolated component of \Spec{A/(I \cap A)}.

It is of interest computationally to be able to replace
multiplicatively closed sets by powers of a single element. In the
proof of \propref{initial:genericprop}, we have chosen an element $r$
which vanishes on every primary component of \Spec{A/(I \cap A)}
except the one defined by $q$. By construction, the product of $r$
with any element of $q$ vanishes everywhere. We observe that our
choice of a single element $r$ differs from the construction of single
elements to replace multiplicatively closed sets in \cite{gtz88},
Prop.\ 3.7:

\begin{example}
\label{initial:gtzex}
Let $A = k[a, b]/(a b)$, let $\Ax = A[x]$, and let $I =
\mbox{$(a(a-1)x,\, x^2)$}$. If $p$ is chosen to be the minimal prime
$(b) \subset A$, then $a \not\in p$ and $a \in \Ann{(b)}$. If we let
$r = a$, and let $B = A_r/p_r$, we have $$ \iin{I} \Bx = \iin{I \Bx} =
((a-1)x,\, x^2).$$ On the other hand, $$ \iin{I} \kpx = \iin{I \kpx} =
(x).$$ Following \cite{gtz88}, if we take $s = a (a-1)$, then $$ I
\Apx \cap \Ax = I A_s[\x] \cap \Ax = (x).$$ $r$ cannot replace $s$ in
this role.
\end{example}

In other words, the extension and contraction of an ideal with respect
to a local ring strips away all but generic behavior along the
corresponding prime, while it is possible to specialize a \std\ basis
in the sense of \propref{initial:genericprop} and still retain some
information about nongeneric behavior along the primary component.

 \propref{initial:genericprop} makes no claims about the relationship
between initial ideals and their specializations to specific primes.
It can happen that no specialization to a prime is well-behaved, as is
illustrated by the following example.

\begin{example}
\label{initial:fuzzex}
We modify example \exref{initial:ex1}. Let $A=k[a,b]/(a^2)$, let
$\Ax=A[x,y]$, and let $I = (ax-y)$. For any prime $p \subset A$,
$\iin{I} \kpx = (y^2)$, but $\iin{I \kpx} = (y)$.

\end{example}

We would like to associate a set of monomials with each fiber of the
family defined by $I$, corresponding over each prime $p$ to the
monomials of \iin{I \kpx}, and then assert that \iin{I} encodes enough
information to determine these sets generically, i.e. along a dense
open subscheme of the base. When the family has a nonreduced base, as
in \exref{initial:fuzzex}, what should we do over the fuzz? One feels
in this example that the fiber monomial ideals are generically $(x)$,
i.e. along the open set away from the subscheme cut out by $(a)$.
Alas, this open set is empty. This same phenomenon can be observed in
studying the ``open nature of flatness'', where for a nonreduced base
scheme, the open set along which a family is flat can be empty. One
could think of such open sets as being supported on the fuzz away from
a proper subscheme.

The following proposition characterizes those primes which are certain
to be well behaved with respect to specialization of a given \std\
basis.

\begin{prop}
\label{initial:whenprop}
Let $I \subset \Ax$ be an ideal, let $p \subset A$ be a prime, and let
$B = A_p/(I \cap A)_p$. If for each monomial $\x^E$, $\iin{I}_E B$ is
either $(0)$ or $(1)$, then $$\iin{I} \kpx = \iin{I \kpx}.$$
\end{prop}

\begin{proof}
If $I \cap A \not\subset p$, then $B$ is the zero ring, and $\iin{I}
\kpx = \iin{I \kpx} = (1)$. Otherwise, using \lemref{initial:quolem},
we can reduce to the case where $I \cap A = (0)$, so $B = A_p$. By
\corref{initial:gtzcor}, we know in any case that $$ \iin{I} \Apx =
\iin{I \Apx}.$$ Let $J = I \Apx$; we need to show that $$ \iin{J \kpx}
\subset \iin{J} \kpx.$$

Given $\x^E \in \iin{J \kpx}$, choose $f \in J \kpx$ so $\iin{f} =
\x^E$. Among all $g \in J$ with image $f$ in $k(p)[\x]$, choose one
with minimal leading term \iin{g}. We claim that $\iin{g} = (1+c)\x^E$
with $c \in p$, so $\x^E \in \iin{J}\, k(p)[\x]$.

Suppose otherwise, that $\iin{g} = c\x^D$ with $\x^D > \x^E$. Then $c
\in p$, and $\iin{I}_D A_p = (1)$. Choose $h \in J$ so $\iin{h} =
\x^D$. Then $g - ch$ also has image $f$ in $k(p)[\x]$, and has a lower
leading term than $g$, contradicting the minimality of our choice for
$g$.
\end{proof}

Let $X \subset \Spec{A}$ be the support of the family defined by $I$;
$X$ is cut out by $I \cap A$. Geometrically, the criterion of
\propref{initial:whenprop} is satisfied if the zero locus of each
coefficient ideal of \iin{I} either avoids the point $p$, or contains
an open neighborhood of $p$ in $X$.

This criterion is sufficient, but not necessary, for taking initial
ideals to commute with specialization. For example, if $I$ is a
monomial ideal, it can have arbitrary coefficient ideals, yet $\iin{I
\kpx} = \iin{I} \kpx$ for all primes $p$.


\section{Faithful Flatness}
\label{flat}
\setcounter{defn0}{0}

The criterion of \propref{initial:whenprop} gives a sufficient
condition for $\Ax/I$ to be faithfully flat over $A/(I \cap A)$. We
will apply the following criterion for faithful flatness; see
Matsumura \cite{mat86}, Atiyah and MacDonald \cite{atma69}, or
Bourbaki \cite{bou89}, for full expositions.

\begin{lemma}
\label{flat:lemma}
Let $M$ be an $A$-module. If for each prime $p \subset A$, $M_p$ is a
nontrivial, free $A_p$-module, then $M$ is faithfully flat over $A$.
\qed
\end{lemma}

The following lemma will be used in two different proofs.

\begin{lemma}
\label{flat:nmlem}
Let $I \subset \Ax$ be an ideal, and define $M = \Ax/I$. Let $V = \{
\x^E \mid \iin{I}_E \neq (1) \}$, and let $N \subset M$ be the
$A$-submodule generated by $V$. Then $N = M$.
\end{lemma}

\begin{proof}
Suppose that $N \neq M$, and choose a nonzero element $f \in M/N$.
Among all $g \in \Ax$ with image $f$ in $M/N$, choose one with minimal
leading term \iin{g}. Let $\iin{g} = c\x^F$. If $\x^F \in V$, then
$c\x^F \in N$, so $g - c\x^F$ also represents $f$, contradicting the
minimality of our choice for $g$. On the other hand, if $\x^F \not\in
V$, then $\iin{I}_F = (1)$, so $\iin{h} = \x^F$ for some $h \in I$.
Then $g - ch$ also represents $f$, again contradicting the minimality
of our choice for $g$.
\end{proof}

\begin{prop}
\label{flat:prop1}
Let $I \subset \Ax$ be a proper ideal, and define $M = \Ax/I$. If for
each prime $p \subset A$ and for each monomial $\x^E$, $\iin{I}_E B$
is either $(0)$ or $(1)$ where $B = A_p/(I \cap A)_p$, then $M$ is a
faithfully flat $A/(I \cap A)$-module.
\end{prop}

\begin{proof}
Disregard primes $p \not\supset I \cap A$. Using
\lemref{initial:quolem}, we can reduce to the case where $I \cap A =
(0)$. We want to show that $M$ is a faithfully flat $A$-module.

Given a prime $p \subset A$, let $V = \{ \x^E \mid \iin{I}_E A_p = (0)
\}$, and let $N \subset M_p$ be the $A_p$-submodule generated by $V$.
$N = M_p$ by \lemref{flat:nmlem}; we claim that $N$ is nontrivial and
free. The result then follows from \lemref{flat:lemma}.

$N$ is nontrivial because $I \cap A = (0)$, so $\iin{I}_1 = (0)$, and
$1 \in V$. $N$ is free, because any relation among its generators
would be an element of $I \Apx$, whose lead term belongs to $V$. Since
$\iin{I \Apx}_E = \iin{I}_E A_p$ by \corref{initial:gtzcor}, this
would contradict the definition of $V$.
\end{proof}

\begin{corollary}
\label{flat:cor}
Given a \std\ basis $\{ f_1,\, \ldots,\, f_r \}$ for $I \subset \Ax$,
let $T$ denote the finite set of exponents $E$ for which $\iin{I}_E
\neq I \cap A$, and which occur as the exponent of some \iin{f_i}. If
$$s \in \bigcap_{E \in T}\; \{ \sqrt{(J^2:J)} \mid J = \iin{I}_E \}$$
and $s \not\in \sqrt{I \cap A}$, then $M_s$ is faithfully flat over
$A_s/(I \cap A)_s$.
\end{corollary}

\begin{proof}
We have the equality of sets
$$ \{ \iin{I}_E \mid \iin{I}_E \neq I \cap A \} \; =
\{ \sum_i \iin{I}_{E_i} \mid E_i \in T \; \mbox{for each} \; i \}. $$
Thus, we get the same intersection of ideals if we replace the index
set $T$ by the infinite set of exponents
$$\{ E \mid \iin{I}_E \neq I\cap A\}.$$
Reduce to the case where $I \cap A = (0)$. The second condition, that
$s$ is not nilpotent, insures that $A_s$ is not the zero ring.

For each $J = \iin{I}_E$, $(J^2:J)$ is supported on precisely those
primes $p$ so $J A_p$ is neither $(0)$ nor $(1)$: $$(J^2:J) A_p = (1)
\Leftrightarrow J A_p = J^2 A_p \Leftrightarrow J A_p = (0)
\mbox{\,or\,} (1).$$
\end{proof}

If the scheme defined by $I \cap A$ has any reduced components, then
such $s$ exist. \corref{flat:cor} remains true for $$s \in \bigcap_{E
\in T}\; \{ \sqrt{J} \mid J = \iin{I}_E \},$$ but nontrivial such $s$
need not exist when $I \cap A$ has more than one reduced component, as
is illustrated by the following example.

\begin{example}
\label{flat:redex}
Let $A = k[a,b,c,d]/((a,b) \cap (c,d))$, let $\Ax = A[x,y]$, let $I =
(ax+b, cy+d)$, and let $M = \Ax/I$. $\{ax+b, cy+d\}$ is a \std\ basis
for $I$, and we have $\iin{I}_x = (a)$, $\iin{I}_y = (c)$, $(a^2:a) =
(a,c,d)$, and $(c^2:c) = (a,b,c)$. For any $s \in (a,c,d) \cap (a,b,c)
= (a,c)$, and any prime $p \subset A$ such that $s \not\in p$,
$\iin{I}_x A_p$ and $\iin{I}_y A_p$ are each either $(0)$ or $(1)$.
Thus $M_s$ is faithfully flat over $A_s$ for any such $s$.

 \Spec{A} consists of the union of the two planes $a = b = 0$ and $c =
d = 0$, meeting at a common origin. $\iin{I}_x = (a)$ vanishes
identically on the plane $a = b = 0$, and is nonzero away from the
line $a = 0$ on the plane $c = d = 0$. Thus $(a^2:a)$ is supported on
this line, which is the locus where $\iin{I}_x$ is locally neither
$(0)$ nor $(1)$. Analogous statements hold for $\iin{I}_y$.

In this example, $$\sqrt{\iin{I}_x} \cap \sqrt{\iin{I}_y} = (a) \cap
(c) = (0),$$ so we cannot simplify the criterion of \corref{flat:cor}.
\end{example}

Compare \propref{flat:prop1} with \propref{initial:whenprop}, and
\exref{initial:ex1}. Summarizing, if the set of monomials of $J =
\iin{I} \kpx$ is locally constant as a function of $p$, then \std\
bases are well behaved with respect to specialization, and moreover,
$M$ is a faithfully flat $A$-module. However, away from the locus
where $J$ is locally constant, $M$ can remain faithfully flat over
$A$.

The next proposition gives some geometric consequences of faithful
flatness. Recall that a morphism of schemes $f: X \rightarrow Y$ is
said to be surjective if for every point $p \in Y$, there exists a
point $P \in X$ such that $f(P) = p$. A morphism of schemes is said to
be dominant if for every point $P \in X$, the induced map $f_P^\#:
\Oh_{Y,f(P)} \rightarrow \Oh_{X,P}$ is injective. While the first
condition is purely topological, the second condition considers the
effect of $f$ on the sheaves of rings $\Oh_X$, $\Oh_Y$, and does not
imply the first.

\begin{prop}
\label{flat:prop2}
Let $I \subset \Ax$ be an ideal, and let $X = \Spec{\Ax/I}$, $Y =
\Spec{A}$. If $M = \Ax/I$ is a faithfully flat $A$-module, then the
corresponding morphism of schemes $X \rightarrow Y$ is surjective and
dominant, and maps each primary component of $X$ to a primary
component of $Y$.
\end{prop}

\begin{proof}
Once a map is known to be flat, the surjectivity of $X \rightarrow Y$
is an equivalent definition of faithful flatness (\cite{mat86}, Thm.\
7.3; \cite{atma69}, Ch.\ 3., Ex.\ 16). Given any prime $P \subset M$,
let $p = P \cap A$. The local homomorphism $A_p \rightarrow M_P$ is
faithfully flat, and thus injective, because $M$ is flat over $A$
(\cite{atma69}, Ch.\ 3., Ex.\ 18). Thus $X \rightarrow Y$ is dominant.
It remains to prove that if $P$ is an associated prime of $(0)$ in
$M$, then $p = P \cap A$ is an associated prime of $(0)$ in $A$. This
follows from \cite{bou89}, Ch.\ 4, \S 2.6, Cor.\ 1 to Thm.\ 2.
\end{proof}

Combining \propref{flat:prop1} with \propref{flat:prop2}, \std\ bases
can be used to manipulate the component structure of the total space
of a family of schemes. In \cite{gtz88}, this problem is approached
differently, via rings of fractions. For example, in Cor.\ 3.8 of
\cite{gtz88}, for $A$ an integral domain with quotient field $K$, the
coefficients of \iin{I} are used to find an $s \in A$ so $I\, A_s[\x]
\cap \Ax$ computes $I\, K[\x] \cap \Ax$; the resulting ideal is then
the intersection of the components of $I$ which surject onto the base
\Spec{A}. In this setting, our choice of $s$ in \corref{flat:cor}
represents a modest improvement over their choice, and
\propref{flat:prop2} illuminates the connection between these
approaches.

When $M$ is finitely generated as an $A$-module, one could determine
the point set in \Spec{A} over which faithful flatness fails, by
studying a presentation matrix of $M$ as an $A$-module; the maximal
minors of this matrix generate one of the {\it Fitting ideals} of $M$.
The use of \std\ bases is more efficient than a brute force study of
these minors, as is evidenced by the following example:

\begin{example}
\label{flat:fitex}
This example is a modification of \exref{initial:ex1}. Let $A=k[a]$,
let $\Ax=A[x,y]$, and let $I = \mbox{$(ax+y,$}\, x^3,\, x^2 y,\, x
y^2,\, y^3) \subset \Ax$.  The coefficient ideals of \iin{I} are given
by the following diagram:

\begin{center}
\begin{tabular}{r|llll}
$y$ & $(1)$ \\
      & $(0)$ & $(1)$ \\
      & $(0)$ & $(a)$ & $(1)$ \\
      & $(0)$ & $(a)$ & $(a)$ & $(1)$ \\
\cline{2-5}
\multicolumn{5}{r}{$x$}
\end{tabular}
\end{center}

As an $A$-module, $M = \Ax/I$ is finitely generated by the set of
monomials having nonunit coefficient ideals, $\{x^2,\, x y,\, y^2,\,
x,\, y,\, 1\}$. These generators have as relations the multiples
$ax+y$,\, $x(ax+y)$, and $y(ax+y)$ of the ideal generator $ax+y$. We
organize this data into the following presentation matrix for $M$:

\begin{center}
\begin{tabular}{rl|ccccccl|}
\multicolumn{2}{l}{} &
$x^2$ & $x y$ & $y^2$ & $x$ & $y$ & $1$ &
\multicolumn{1}{l}{} \\
\cline{3-9}
$a x+y$ & & $0$ & $0$ & $0$ & $a$ & $1$ & $0$ & \\
$x(a x+y)$ & & $a$ & $1$ & $0$ & $0$ & $0$ & $0$ & \\
$y(a x+y)$ & & $0$ & $a$ & $1$ & $0$ & $0$ & $0$ & \\
\cline{3-9}
\end{tabular}
\end{center}

The $(x y,\, y^2,\, y)$-minor of this matrix is nonsingular, so $M$ is
flat over $A$. However, the leading nonzero minor, on columns $(x^2,\,
xy,\, x)$, has determinant $a^3$. This minor demonstrates that $M_a$
is flat over $A_a$, but leaves open the question of what happens over
$a = 0$.

The coefficient ideals of \iin{I} determine a locus away from which
this leading minor is nonsingular. Since one only needs to consider
coefficient ideals corresponding to minimal generators of \iin{I},
\std\ bases can be used to find this locus without explicitly
considering every row of the presentation matrix: a single element of
the \std\ basis stands in for many rows of the presentation matrix.

Note also that the coefficient ideals and the determinant give
different scheme structures for the set where this leading minor loses
rank. The coefficient ideals describe the support of the module
defined by the leading minor, while the determinant describes a
thicker scheme enjoying a universal property with respect to base
change; see \cite{eis89}, Ch.\ 10.
\end{example}


\section{Fibers}
\label{fibers}
\setcounter{defn0}{0}

The following pair of propositions concern the behavior of coefficient
ideals with respect to the geometry of the fibers of a family.

\begin{prop}
\label{fibers:isoprop}
Let $I \subset \Ax$ be an ideal, let $M = \Ax/I$, and let $p \subset
A$ be a prime ideal. The following two statements are equivalent:

(a) $f: A_p \rightarrow M_p$ is surjective;

(b) $\iin{I}_{x_i} A_p = (1)$ for each $i$.
\end{prop}

\begin{proof}
If $\iin{I}_{x_i} A_p = (1)$ for each $i$, then $M_p$ admits a
relation of the form $x_i - f_i(x_{i+1},\ldots,x_n)$, for each $i$.
This proves that $f$ is surjective. Conversely, if $f$ is surjective,
then $M_p$ admits a relation of the form $x_i - c_i$ with $c_i \in
A_p$, for each $i$. Thus, the corresponding coefficient ideals are
unit ideals.
\end{proof}

\begin{corollary}
\label{fibers:cor}
Let $I \subset \Ax$ be an ideal, let $X = \Spec{A/(I \cap A)}$, and
let $Y = \Spec{\Ax/I}$. If $I \cap A$ is a prime ideal, and each
coefficient ideal $\iin{I}_{x_i} \neq I \cap A$, then the induced
morphism of schemes $Y \rightarrow X$ is an isomorphism over a
nonempty open subset of $X$.
\qed
\end{corollary}

Let $\iin{I}_{x_i^\infty}$ denote the stationary limit of the
ascending chain of coefficient ideals $\iin{I}_{x_i} \subset
\iin{I}_{x_i^2} \subset \ldots\,$. Each $\iin{I}_{x_i^\infty}$ can be
computed as the ideal generated by the leading coefficients of all
\std\ basis elements having a leading monomial of the form $x_i^e$.

\begin{prop}
\label{fibers:finprop}
Let $I \subset \Ax$ be an ideal, let $M = \Ax/I$, and let $p \in A$ be
a prime ideal. The following two statements are equivalent:

(a) $f: A_p \rightarrow M_p$ is a finite map;

(b) $\iin{I}_{x_i^\infty} A_p = (1)$ for each $i$.
\end{prop}

\begin{proof}
Suppose that (a) holds, so $M_p$ is a finite $A_p$-module. Fix $i$,
and let $N \subset M_p$ be the subalgebra generated by $x_i$. Then $N$
can be generated as an $A_p$-module by the finite set $\{ 1, x_i,
\ldots ,x_i^{e-1}\}$ for some $e$. $x_i^e$ can be expressed in terms
of these generators, yielding an expression in $I \Apx$ with leading
term $x_i^e$. Thus $\iin{I \Apx}_{x_i^e} = (1)$.

Conversely, assume (b). For each $i$, choose $e_i$ so $\iin{I
\Apx}_{x_i^e} = (1)$, and let $V$ denote the finite set of monomials
which do not belong to the ideal $(x_1^{e_1},\ldots,x_n^{e_n})$. If $N
\subset M_p$ is the $A_p$-submodule generated by $V$, then $N$ is
finitely generated, and $N = M_p$ by \lemref{flat:nmlem}.
\end{proof}

Geometrically, let $X = \Spec{\Ax/I}$, let $Y = \Spec{A/(I \cap A)}$,
and let $g: X \rightarrow Y$ be the morphism of schemes induced by the
inclusion $A/(I \cap A) \subset \Ax/I$. If we let $U \subset Y$ be the
complement of the union of the subschemes of $Y$ cut out by the
coefficient ideals $\iin{I}_{x_i}$, then \propref{fibers:isoprop}
asserts that $U$ is the largest open set with the property that
$g^{-1}(U) \rightarrow U$ is an isomorphism. If instead we let $U
\subset Y$ be the complement of the union of the subschemes of $Y$ cut
out by the coefficient ideals $\iin{I}_{x_i^\infty}$, then
\propref{fibers:finprop} asserts that $U$ is the largest open set with
the property that $g^{-1}(U) \rightarrow U$ is a finite morphism.
Thus, while \std\ bases can only generically detect faithful flatness,
they are capable of detecting precisely the locus where a morphism
restricts to an isomorphism, or to a finite map.

Detection of quasi-finite morphisms is more subtle; a family which
restricts to a finite family over an open subset of the base need not
be quasi-finite:

\begin{example}
\label{fibers:ex1}
Let $A=k[a, b]$, let $\Ax=A[x]$, let $I = (a x - b)$, and let $M =
\Ax/I$. The localization $M_a$ is finite over $A_a$. However, the
total space $a x - b = 0$ is irreducible, and consists of a line over
$a = b = 0$.
\end{example}





\begin{thebibliography}{[BaSt 86a]}

\bibitem[AtMa 69]{atma69}
	M. F. Atiyah and I. G. MacDonald,
	{\em Introduction to Commutative Algebra},
	Addison-Wesley Series in Mathematics, 1969.
	ISBN 0-201-00361-9.

\bibitem[Bou 89]{bou89}
	N. Bourbaki,
	{\em Commutative Algebra, Chapters 1-7},
	Springer-Verlag, 1989.
	ISBN 0-387-19371-5.

\bibitem[Eis 89]{eis89}
	D. Eisenbud,
	{\em Commutative Algebra with a view toward Algebraic Geometry},
	Brandeis University Lecture notes, March 1, 1989.

\bibitem[Gal 79]{gal79}
	A. Galligo,
	{\em Theoreme de division et stabilite en geometrie analytique locale},
	Ann. Inst. Fourier  T. 29,2, 1979, pp. 107-184 .

\bibitem[GTZ 88]{gtz88}
	P. Gianni, B. Trager, G. Zacharias,
	Gr\"{o}bner bases and primary decomposition of polynomial ideals,
	{\em J.\,Symbolic Computation\,} {\bf 6}, 1988, pp. 148-166,
	a.k.a. {\em Computational Aspects of Commutative Algebra},
	ed. L. Robbiano, Academic Press, San Diego, 1989, pp. 15-33.
	ISBN 0-12-589590-9.

\bibitem[HLT 73]{hlt73}
	H. Hironaka, M. Lejeune-Jalabert, B. Teissier,
	{\em Platificateur local en geometrie analytique et applatissement local},
	Asterisque n. 7 \& 8, ``Singularites a Cargese'', 1973, pp. 441-466.

\bibitem[Mat 86]{mat86}
	H. Matsumura,
	{\em Commutative ring theory},
	Cambridge University Press, Cambridge, 1986.
	ISBN 0-521-36764-6.

\end{thebibliography}
\end {document}